\documentclass{amsart}
\def\eps{\varepsilon}
\newtheorem{lemma}{Lemma}[section]
\newtheorem{theorem}[lemma]{Theorem}
\newtheorem{remark}[lemma]{Remark}

\newtheorem{coro}[lemma]{Corollary}
\newtheorem{definition}[lemma]{Definition}
\newtheorem{example}[lemma]{Example}
\newcommand{\been}{\begin{enumerate}}
\newcommand{\enen}{\end{enumerate}}
\newcommand{\beeq}{\begin{eqnarray}}
\newcommand{\eneq}{\end{eqnarray}}
\newcommand{\beeqn}{\begin{eqnarray*}}
\newcommand{\eneqn}{\end{eqnarray*}}
\newcommand{\bd}{\begin{displaymath}}
\newcommand{\ed}{\end{displaymath}}

\title[Generalization of the Second Bogolyubov's Theorem]
{Generalization of the Second Bogolyubov's Theorem for
Non-Almost Periodic Systems}

\author{David N.~Cheban}
\address[D. Cheban]{%
State University of Moldova\\ Department of Mathematics and
Informtics\\ A. Mateevich Street 60\\ MD--2009 Chi\c{s}in\u{a}u,
Moldova} \email[D. Cheban]{cheban@usm.md}
\author{Jinqiao  Duan}
\address[J. Duan]{%
Department of Applied Mathematics \\ Illinois Institute of
Technology\\ Chicago, IL 60616, USA} \email[J.~Duan]{duan@iit.edu}
\author{Anatoly Gherco}
\address[A. Gherco]{%
State University of Moldova\\ Department of Mathematics and
Informtics\\ A. Mateevich Street 60\\ MD--2009 Chi\c{s}in\u{a}u,
Moldova} \email[A. Gherco]{gerko@usm.md}
\date{\today}
\subjclass{primary:34D20, 34D40, 34D45, 58F10, 58F12, 58F39;
secondary: 35B35, 35B40.} \keywords{Nonautonomous dynamical
system, skew--product flow, recurrent and almost periodic solutions, uniform
averaging principle}

\begin{document}
\begin{abstract}
{\bf Nonlinear Analysis B,  4(2003), 599-613.}

 {  The article is  devoted to the generalization of the
second Bogolyubov's theorem to non-almost periodic dynamical
systems. We prove the analog of the
 second  Bogolyubov's theorem for
  recurrent or  pseudo recurrent  dynamical systems  in Banach spaces.
Namely, we obtain the relation between a recurrent dynamical
system and its averaged  dynamical system. We also  study existence
of recurrent and pseudo recurrent motions (including special cases of periodic,
quasi-periodic and almost periodic motions) in related
nonautonomous systems.}
\end{abstract}
\maketitle

\section{Introduction}

The problem of averaging  in time is well-studied for almost
periodic systems in Banach spaces. A well-known result in
this direction is the
  second  Bogolyubov's theorem (see for example \cite{Bog,DK}) which affirms
  that the equation
\begin{equation}\label{eq0.1}
\dot{x}=\varepsilon f(t,x)
\end{equation}
with almost periodic function $f$ for sufficient small
$\varepsilon$ admits a unique  almost periodic solution in the
neighborhood of hyperbolic stationary point $x_0$ of the ``averaged"
equation
\begin{equation} \label{eq0.2}
\dot{x}=\varepsilon f_0(x),
\end{equation}
where
\begin{equation} \label{eq0.3}
 f_0(x)=\lim\limits_{T\to+\infty}\frac{1}{T}\int\limits_t^{t+T}f(s,x)ds
\end{equation}
and the limit (\ref{eq0.3}) is uniformly with respect to  (w.r.t.) $t\in \mathbb R$.
The first Bogolyubov's theorem determines the closeness or nearness
of the solutions on finite time intervals for the original equation
(\ref {eq0.1}) and the
averaged  equation   (\ref {eq0.2}).
Note that periodic and quasi-periodic functions are special almost periodic
functions.

In this paper, we generalize the second  Bogolyubov's theorem
   for the equation (\ref{eq0.1}) to the case when
function $f$ is recurrent or pseudo recurrent (see definitions in Sections
4 and 6).

The paper is organized as follows. In section 2 we study the
 existence of invariant  integral
manifolds  of quasilinear nonautonomous dynamical systems (Theorem
 \ref{t2.5} and Theorem  \ref{t2.6}).
Results in this section are used in the following sections.

Section 3 contains the main results about generalization of the second
Bogolyubov's theorem for non-almost periodic systems (Theorem
\ref{t3.3}  and Theorem \ref{t3.4}).

In section 4 we give conditions of existence of recurrent
solutions of nonautonomous equations in a standard form, if
corresponding averaging equation admits a hyperbolic stationary
point (Theorem \ref{t4.6}).

Second 5 is devoted to study of the  existence of
invariant torus and quasi-periodic solutions
of quasilinear equations on the torus (Theorem \ref{t5.2}, Theorem \ref{t5.3}
and Corollary \ref{c5.4}).

In second 6 we discuss the  existence of pseudo
recurrent integral manifolds (Theorem \ref{soudo}).

\section{Quasilinear nonautonomous dynamical systems }

Let $\Omega$ be a compact metric space and
$(\Omega,{\mathbb R}, \sigma)$ be an autonomous
dynamical system on $\Omega$. Let $E$ be a Banach space,
and $Y$ and $W$ are two complete metric spaces.
Denote $L(E)$  the space of
all   linear continuous operators on $E$ and $C(Y,W)$   the space of
all   continuous functions ${f}:{Y}\to{W}$ endowed by compact-open
topology, i.e., uniform convergence on compact subsets in  $Y$.
We use these notations for the rest of the paper.
The results in this section will be used in later sections.

Consider the linear equation
\begin{equation} \label{eq2.1}
\dot{x}=A(\omega t)x \quad (\omega\in\Omega,\,\, \omega t=\sigma(t,\omega))
\end{equation}
and the inhomogeneous equation
\begin{equation} \label{eq2.2}
\dot{x}=A(\omega t)x+f(\omega t),
\end{equation}
where  $A\in C(\Omega, L(E))$ and   $f\in C(\Omega, E)$.

\begin{definition}\label{d2.1}
Let $U(t,\omega)$ be the operator of Cauchy (solution operator)
of the linear equation $(\ref{eq2.1})$.
The equation (\ref{eq2.1}) is called hyperbolic if there exist positive
numbers
$N$, $\nu>0$ and continuous projection $P\in C(\Omega,L(E))$ $($i.e.
$P^2(\omega)=P(\omega)$ for all $\omega\in \Omega)$ such that
\begin{enumerate}
\item[$1)$] For all $t\in\mathbb R$ and $\omega\in\Omega$,
$U(t,\omega)P(\omega)=P(\omega t)U(t,\omega)$;
\item[$2)$] For all $t\ge \tau$ and $\omega\in\Omega$,
$\Vert U(t,\omega)P(\omega)U^{-1}(\tau,\omega)\Vert\le N\exp{(-\nu(t-\tau))}$;
\item[$3)$] For all $t\le \tau$ and $\omega\in\Omega$,
$\Vert U(t,\omega)Q(\omega)U^{-1}(\tau,\omega)\Vert\le N\exp{(\nu(t-\tau))}$,
 where $Q(\omega)=I-P(\omega))$.
\end{enumerate}
\end{definition}

\begin{definition}\label{d2.2}
The function $G:{\mathbb R}_*^2\times \Omega\to L(E)$ defined by
\begin{equation}\label{eq2.3}
G(t,\tau,\omega)=\left\{\begin{array}{rl}
U(t,\omega)P(\omega)U^{-1}(\tau,\omega)&\,\, \mbox{for}\,\, t>\tau\\
-U(t,\omega)Q(\omega)U^{-1}(\tau,\omega)&\,\, \mbox{for}\,\, t<\tau
\end{array}
\right.
\end{equation}
is called the Green's function for hyperbolic linear equation $(\ref{eq2.1})$,
where ${\mathbb R}_*^2={\mathbb R}^2\setminus {\Delta}_{{\mathbb R}^2}$,
${\Delta}_{{\mathbb R}^2}=\{(t,t)\mid\, t\in {\mathbb R}\}$ and
$P$, $Q$ are the projections from definition $\ref{d2.1}$.
\end{definition}

\begin{remark}\label{r2.2}
The Green's function satisfies the following conditions
$($see \cite{DK} and \cite{Sam}$)$:
\begin{enumerate}
\item[$1)$] For every $t\ne\tau$ the function $G(t,\tau,\omega)$ is
continuously differentiable and
\begin{displaymath}
\frac{\partial {G(t,\tau,\omega)}}{\partial {t}}=A(\omega t)G(t,\tau,\omega)
\quad (\omega\in\Omega).
\end{displaymath}
\item[$2)$] $G(\tau+0,\tau,\omega)-G(\tau-0,\tau,\omega)=I$
$(\tau\in\mathbb R$, $\omega\in\Omega)$.
\item[$3)$] $\Vert G(t,\tau,\omega)\Vert\le N\exp{(-\nu\vert t-\tau\vert)}$
$(t,\tau\in\mathbb R$, $\omega\in\Omega)$.
\item[$4)$] $G(0,\tau,\omega t)=G(t,t+\tau,\omega)$
$(t,\tau\in\mathbb R$, $\tau\ne0$, $\omega\in\Omega)$.
\end{enumerate}

\end{remark}
\begin{theorem}\label{t2.4}
Suppose that the   linear equation $(\ref{eq2.1})$ is hyperbolic.
Then for
$f\in C(\Omega, E)$, the function $\gamma(\omega)$ defined by
\begin{equation} \label{eq2.4}
\gamma(\omega)=\int\limits_{-\infty}^{+\infty}G(0,\tau,\omega)f(\omega \tau)d\tau
\quad (\omega\in\Omega)
\end{equation}
is   continuous, i.e.,   $\gamma\in C(\Omega, E)$,  and
\begin{enumerate}
\item[$1)$]  $\gamma(\omega t)= \varphi(t,\gamma(\omega),\omega)$
holds for all $\omega\in\Omega$ and
$t\in\mathbb R_{+}$, where $\varphi(t,x,\omega)$ is the unique solution of
the corresponding inhomogeneous equation $(\ref{eq2.2})$ with the initial condition $\varphi(0,x,\omega)=x$;
\item[$2)$] $\Vert\gamma\Vert\le\frac{2N}{\nu}\Vert f \Vert$, where
$\Vert\gamma\Vert=\max\limits_{\omega\in\Omega}\vert\gamma(\omega)\vert$.
\end{enumerate}
\end{theorem}
\begin{proof}
The proof of this assertion is obtained by slight modification of
arguments from \cite[Ch.III]{Sam} and we omit the details.
\end{proof}

\vspace{2mm}

Let us consider the following  quasilinear equation in Banach space $E$
\begin{equation} \label{eq2.5}
\dot{x}=A(\omega t)x+f(\omega t)+F(\omega t,x),
\end{equation}
where  $A\in C(\Omega, L(E))$,   $f\in C(\Omega, E)$ and $F\in C(\Omega\times E, E)$.

\begin{theorem}\label{t2.5}
(Invariant integral manifold)
Assume that there exist positive numbers $L<L_0 := \frac{\nu}{2N}$
and $r<r_0:=\gamma_0(\frac{\nu}{2N}-L_0)^{-1}$ such that
\begin{equation} \label{eq2.6}
\Vert F(\omega,x_1)-F(\omega,x_2) \Vert\le L\Vert x_1-x_2\Vert
\end{equation}
for all $\omega\in\Omega$ and
$x_1,x_2\in B[Q,r]=\{x\in E\mid \rho(x,Q)\le r\}$,
where $Q=\gamma(\Omega)$, $\gamma\in C(\Omega,E)$ is defined in
$(\ref{eq2.4})$ and $\gamma_0=
\max\limits_{\omega\in\Omega}\Vert F(\omega,\gamma(\omega))\Vert$.
Then there exists a unique function $u\in C(\Omega,B[Q,r])$ such that
\begin{equation} \label{eq2.7}
u(\omega t)=\psi(t,u(\omega),\omega)
\end{equation}
for all $t\in\mathbb R_{+}$ and $\omega\in\Omega $, where
$\psi(\cdot,x,\omega)$
is the unique solution of the quasilinear
equation $(\ref{eq2.5})$ with the initial condition
$\psi(0,x,\omega)=x$. Therefore, the graph of $u$ is
an invariant integral manifold for the quasilinear
equation $(\ref{eq2.5})$.
\end{theorem}

\begin{proof}
Let $x=y+\gamma(\omega t)$. Then from the equation $(\ref{eq2.5})$ we obtain
\begin{equation} \label{eq2.8}
\dot{y}=A(\omega t)y+F(\omega t,y+\gamma(\omega t)).
\end{equation}

If $0<r<r_0$ and $\alpha\in C(\Omega,B[Q,r])$,  then   the equality
\begin{equation} \label{eq2.9}
(\Phi\alpha)(\omega)=\int\limits_{-\infty}^{+\infty}
G(0,\tau,\omega)F(\omega \tau,\alpha(\omega\tau)+\gamma(\omega\tau))d\tau
\end{equation}
defines a function $\Phi\alpha\in C(\Omega,E)$.
In virtue of Theorem \ref{t2.4},  we have
\begin{eqnarray}
\label{eq2.10}
&{}&\Vert\Phi\alpha\Vert\le \frac{2N}{\nu}\max\limits_{\omega\in\Omega}
\Vert F(\omega,\alpha(\omega)+\gamma(\omega))\Vert\notag\\ &{}&
\le \frac{2N}{\nu}\max\limits_{\omega\in\Omega}
\Vert F(\omega,\alpha(\omega)+\gamma(\omega))-
F(\omega,\gamma(\omega))\Vert + \frac{2N}{\nu}\max\limits_{\omega\in\Omega}
\Vert F(\omega,\gamma(\omega))\Vert\notag\\ &{}&
\le\frac{2N}{\nu}L\Vert \alpha\Vert+ \frac{2N}{\nu}\gamma_0\le
\frac{2N}{\nu}Lr+ \frac{2N}{\nu}\gamma_0\le
\frac{2N}{\nu}L_0r_0+ \frac{2N}{\nu}\gamma_0=r_0
\end{eqnarray}

and consequently $\Phi(C(\Omega,B[Q,r_0]))\subseteq C(\Omega,B[Q,r_0])$.

Now we will show that the mapping $\Phi:C(\Omega,B[Q,r_0])\to
C(\Omega,B[Q,r_0])$ is Lipschitzian. In fact, according to Theorem
\ref{t2.4} we have
\begin{eqnarray}
\label{eq2.11}
&{}&\Vert\Phi\alpha_1-\Phi\alpha_2\Vert\le \frac{2N}{\nu}\max\limits_{\omega\in\Omega}
\Vert F(\omega,\alpha_1(\omega)+\gamma(\omega))-
F(\omega,\alpha_2(\omega)+\gamma(\omega))\Vert\notag\\&{}&
\le \frac{2N}{\nu}L\max\limits_{\omega\in\Omega}
\Vert \alpha_1(\omega)-\alpha_2(\omega)\Vert.
\end{eqnarray}

We note that $\frac{2N}{\nu}L\le\frac{2N}{\nu}L_0<1$. Thus the mapping
$\Phi$ is a contraction and, consequently by Banach fixed point theorem,
 there exists a unique function
$\alpha\in C(\Omega,B[Q,r_0])$ such that $\Phi\alpha=\alpha$. To finish the
proof of the theorem it is sufficient to put $u=\gamma+\alpha$.
\end{proof}

We now consider the perturbed quasilinear equation
\begin{equation} \label{eq2.12}
\dot{x}=A(\omega t)x+f(\omega t)+\varepsilon F(\omega t,x),
\end{equation}
where $\varepsilon\in[-\varepsilon_0,\varepsilon_0]$ $(\varepsilon_0>0)$
is a small parameter. We have a similar theorem.

\begin{theorem}\label{t2.6}
(Invariant integral manifold and convergence)
Assume that there exist positive numbers $r$ and $L$ such that
\begin{equation} \label{eq2.13}
\Vert F(\omega,x_1)-F(\omega,x_2) \Vert\le L\Vert x_1-x_2\Vert
\end{equation}
for all $\omega\in\Omega$ and $x_1,x_2\in B[Q,r]$.
Then for sufficiently small $\varepsilon$
there exists a unique function $u_{\varepsilon}\in C(\Omega,B[Q,r])$ such that
\begin{equation} \label{eq2.14}
u_{\varepsilon}(\omega t)=\psi_{\varepsilon}(t,u_{\varepsilon}(\omega),\omega)
\end{equation}
for all $t\in\mathbb R_{+}$ and $\omega\in\Omega$, where
$\psi_{\varepsilon}(\cdot,x,\omega)$
is the unique solution of equation $(\ref{eq2.12})$ with the initial condition
$\psi_{\varepsilon}(0,x,\omega)=x$.
Moreover,
\begin{equation} \label{eq2.15}
\lim\limits_{\varepsilon\to 0}\max\limits_{\omega\in\Omega}
\Vert u_{\varepsilon}(\omega)-\gamma(\omega)\Vert=0,
\end{equation}
where $\gamma\in C(\omega,E)$ is defined in  $(\ref{eq2.4})$.
\end{theorem}

\begin{proof} We can prove the existence of $ u_\eps $ by slight
modification of the proof of Theorem \ref{t2.4}.

To prove $(\ref{eq2.15})$ we note that
\begin{equation} \label{eq2.16}
\Vert F(\omega,u_{\varepsilon}(\omega)\Vert\le \Vert
F(\omega,u_{\varepsilon}(\omega)-F(\omega,\gamma(\omega)\Vert
+\Vert F(\omega,\gamma(\omega)\Vert\le L r+\gamma_0
\end{equation}
and
\begin{equation} \label{eq2.17}
\Vert u_{\varepsilon}(\omega)-\gamma(\omega)\Vert\le
\Vert \int\limits_{-\infty}^{+\infty}\varepsilon G(0,\tau,\omega \tau)
F(\omega \tau,u_{\varepsilon}(\omega \tau))d\tau
\le\vert \varepsilon\vert \frac{2N}{\nu}(L r+\gamma_0)
\end{equation} 
for all $\omega\in\Omega$ and $\varepsilon\in[-\varepsilon_0,\varepsilon_0]$.
Passing to the limit in the inequality (\ref{eq2.17}) as $\varepsilon\to 0$
we obtain (\ref{eq2.15}).
\end{proof}

\section{Generalization of second Bogolyubov's theorem
for non-almost periodic systems }

In this section, we consider an analog of the second Bogolyubov's theorem
 for the nonautonomous system
\begin{equation} \label{eq3.1}
\dot{x}=\varepsilon f(\omega t,x),
\end{equation}
where $\varepsilon\in[0,\varepsilon_0]$ $(\varepsilon_0>0)$ is a small
parameter.  We do {\em not} assume that  $f$ is almost periodic in  time $t$.
Suppose that the   averaging
\begin{equation} \label{eq3.2}
 \overline{f}(x)=\lim\limits_{T\to+\infty}\frac{1}{T}\int\limits_0^{T}f(\omega t,x)dt
\end{equation}
 exists uniformly w.r.t.   $\omega\in\Omega$,
 and also  uniformly w.r.t. $x$ on every  bounded subset of $E$.

 \begin{remark}\label{r3.1}
 The condition $(\ref{eq3.2})$ is fulfilled if a dynamical system
 $(\Omega,\mathbb R,\sigma)$ is strictly ergodic, i.e. on $\Omega$
 exists a unique invariant measure $\mu$ w.r.t. $(\Omega,\mathbb R,\sigma)$.
 \end{remark}

Along with equation (\ref{eq3.1}) we consider the averaged equation
\begin{equation} \label{eq3.3}
 \dot x=\varepsilon\overline{f}(x).
\end{equation}
Setting slow time $\tau=\varepsilon t$ $(\varepsilon>0)$,
 the equations (\ref{eq3.1})
and (\ref{eq3.3}) can be written in the following form:
\begin{equation} \label{eq3.4}
 \frac{d x}{d\tau}=f(\omega\frac{\tau}{\varepsilon},x)
\end{equation}
and
\begin{equation} \label{eq3.5}
 \frac{d x}{d\tau}=\overline{f}(x)
\end{equation}
respectively.

Suppose that for certain point $x_0\in E$
\begin{equation} \label{eq3.6}
 \overline{f}(x_0)=0,
\end{equation}
then the equation (\ref{eq3.3}) admits a stationary solution
$\varphi_{\varepsilon}(t,x_0)\equiv x_0$.

Assume that the following conditions are fulfilled:

(i) Function $f\in C(\Omega\times B[x_0,r],E)$ where $B[x_0,r]=\{x\in E\mid
\Vert x-x_0\Vert\le r\}$  and $r>0$, and $F$ is bounded on
$\Omega\times B[x_0,r]$. The limit (\ref{eq3.2}) is uniform w.r.t.
$(\omega,x)\in \Omega\times B[x_0,r]$ and functions
$f'_x(\omega,x)$ and $\overline{f}'(x)$ are bounded on
$\Omega\times B[x_0,r]$.

(ii) Functions $f(\omega,x)$ and $\overline{f}(x)$ are twice
continuously differentiable w.r.t variable $x\in B[x_0,r]$.

(iii) The equality (\ref{eq3.2}) can be twice differentiated, i.e., the
following equalities
\begin{equation} \label{eq3.7}
 \overline{f}^{(k)}(x)=
 \lim\limits_{T\to +\infty}\frac{1}{T}\int\limits_0^Tf_x^{(k)}(\omega t,x)d t
 \quad (k=1,2)
\end{equation}
 hold uniformly w.r.t $\omega\in\Omega$ and $x\in B[x_0,r]$.

 We note that
\begin{equation} \label{eq3.8}
 \overline{f}(x+h)-\overline{f}(x)=\overline{f}'(x)h +R(x,h) \quad (x,x+h\in B[x_0,r])
\end{equation}
 where $\Vert R(x,h)\Vert=o(\Vert h\Vert)$.

 Let $A=\overline{f}'(x_0)$  and $B(h)=R(x_0,h)$. Then according to (\ref{eq3.6})
 and (\ref{eq3.8}) we have
\begin{equation} \label{eq3.9}
 \overline{f}(x+h)=Ah +B(h)
\end{equation}

 It is clear (see \cite[Ch.7]{DK})  that the function $B(h)$ satisfies the
 condition of Lipschitz
\begin{equation} \label{eq3.10}
\Vert B(h_1)-B(h_2)\Vert\le L(r)\Vert h_1-h_2\Vert
\end{equation}
$(h_1,h_2\in B[x_0,r])$ and $L(r)\to 0$ as $r\to 0$.

The equation (\ref{eq3.1}) can be rewritten in the following form
\begin{equation} \label{eq3.11}
\frac{dh}{dt}=\varepsilon Ah+\varepsilon g(\omega t,h),
\end{equation}
where $h=x-x_0$ and
\begin{equation} \label{eq3.12}
g(\omega ,h)=f(\omega,x+h)-\overline{f}(x_0+h)+B(h).
\end{equation}

In the equation (\ref{eq3.12}) we make the following change of variable
\begin{equation} \label{eq3.13}
h=z-\varepsilon v(\omega,z,\varepsilon),
\end{equation}
where
\begin{equation} \label{eq3.14}
v(\omega,z,\varepsilon)=\int\limits_0^{+\infty}V(\omega s,z)\exp(-\varepsilon s)ds
\end{equation}
and
\begin{equation} \label{eq3.15}
V(\omega,z)=f(\omega,x_0+z)-\overline{f}(x_0+z).
\end{equation}

\begin{lemma}\label{l3.1}
(\cite[p.457]{DK})
Let $\varphi:{\mathbb R}_+\times \Lambda\to E$ be a function satisfying the
following conditions:

$1$. $M:=\sup\{\Vert \frac{1}{t}\int\limits_0^{t}\varphi(s,\lambda)ds\Vert\mid
t\ge0,\, \lambda\in\Lambda\}<+\infty$.

$2$. $\lim\limits_{T\to+\infty}\frac{1}{T}\int\limits_0^{T}\varphi(s,\lambda)ds=0$
uniformly w.r.t. variable $\lambda\in\Lambda$.

Then the following equality
$$ \lim\limits_{p\to 0}p\int\limits_0^{+\infty}\varphi(s,\lambda)\exp(-ps)ds=0
$$
takes place uniformly w.r.t. $\lambda\in\Lambda$.
\end{lemma}

\begin{lemma}\label{l3.2}
The following equalities
\begin{equation} \label{eq3.16}
 \lim\limits_{\varepsilon\downarrow 0}\varepsilon v(\omega,z,\varepsilon)=0
\end{equation}
and
\begin{equation} \label{eq3.17}
 \lim\limits_{\varepsilon\downarrow 0}\varepsilon v'_z(\omega,z,\varepsilon)=0
\end{equation}
are fulfilled uniformly w.r.t. $\omega\in\Omega$ and $z\in B[0,r]$.
\end{lemma}

\begin{proof}
This assertion follows from Lemma \ref{l3.1}.  In fact, in virtue of
(\ref{eq3.2}), (\ref{eq3.7}) and (\ref{eq3.14}),  the   bounded
function $V(\omega s, z)$ and $V'_z(\omega s, z)$ satisfy the
conditions of Lemma \ref{l3.1}.
\end{proof}
\vspace{2mm}

From the equality (\ref{eq3.7})  it follows that for sufficiently small $\varepsilon>0$
the operator $I-\varepsilon v'_z(\omega,z,\varepsilon)$
$(\omega\in\Omega$, $z\in B[0,r])$ is invertible and
$(I-\varepsilon v'_z(\omega,z,\varepsilon))^{-1}$ is bounded, and, consequently,
the mapping (\ref{eq3.13}) is invertible. According to   equation (\ref{eq3.16})
in the sufficiently small neighborhood of zero and  for sufficiently small $\varepsilon>0$,
we can make the change of variable (\ref{eq3.13}).

Note that
$$
v(\omega t,z,\varepsilon)=\int\limits_0^{+\infty}V(\omega (t+s),z,\varepsilon))
\exp(-\varepsilon s)ds =
\exp(\varepsilon t)\int\limits_t^{+\infty}V(\omega s,z,\varepsilon))
\exp(-\varepsilon s)ds
$$
and we find that
\begin{equation} \label{eq3.18}
\frac{d}{dt}v(\omega t,z,\varepsilon)=\varepsilon v(\omega t,z,\varepsilon)
-V(\omega t,z,\varepsilon)
\end{equation}
and, consequently,
\begin{equation} \label{eq3.19}
\frac{dh}{dt}=\frac{dz}{dt}-\varepsilon v'_z\frac{dz}{dt}-\varepsilon^2 v+
\varepsilon V.
\end{equation}

Using the relation (\ref{eq3.13}), (\ref{eq3.15}) and (\ref{eq3.19})
we reduce the equation  (\ref{eq3.11}) to the form
\begin{eqnarray}
\label{eq3.20}
&{}& (I-\varepsilon v'_z(\omega t,z,\varepsilon))\frac{dz}{dt}=
\varepsilon[f(\omega t,x_0+h,\varepsilon)-f(\omega t,x_0+z,\varepsilon)]\notag\\ &{}&
+\varepsilon \overline{f}(x_0+z)+\varepsilon^2v(\omega t,z,\varepsilon)=
\varepsilon(Az+B(z))\notag\\ &{}&+
\varepsilon^2v(\omega t,z,\varepsilon) +\varepsilon[f(\omega t,x_0+z-\varepsilon v,
\varepsilon)-f(\omega t,x_0+z,\varepsilon)].
\end{eqnarray}

After multiplication of the both sides of the equation  (\ref{eq3.20}) by
$(I-\varepsilon v'_z(\omega t,z,\varepsilon))^{-1}$ and introduction of the
"slow" time $\tau=\varepsilon t$ we obtain
\begin{equation} \label{eq3.21}
\frac{dz}{d\tau}=Az+F(\omega \frac{\tau}{\varepsilon},z,\varepsilon)
\end{equation}
where $F$ possesses the following properties:

a) $F$ admits  a bounded derivable $F'_z(\omega,z,\varepsilon)$
$(\omega\in\Omega$, $z\in B[0,r]$ and $\varepsilon\in[o,\varepsilon_0])$;

b) $F(\omega,z,\varepsilon)=B(z)+O(z)$
uniformly w.r.t. $\omega\in\Omega$ and $z\in B[0,r]$;

c) For every $M>0$ and $\mu>0$, there exists positive numbers
$\varepsilon'_0\le\varepsilon_0$ and $\beta_0$ such that for
$0<\varepsilon<\varepsilon'_0$, and $\Vert z\Vert<\beta_0$,  the
inequalities
\begin{equation} \label{eq3.22}
\Vert F(\omega,z,\varepsilon)\Vert\le M
\end{equation}
and
\begin{equation} \label{eq3.23}
\Vert F(\omega,z_1,\varepsilon)-F(\omega,z_2,\varepsilon)\Vert
\le \mu\Vert z_1-z_2\Vert
\end{equation}
take place for all
$\omega\in\Omega$, $z_1,z_2\in B[0,r]$ and $0\le\varepsilon\le\varepsilon'_0$.

\begin{theorem}\label{t3.3}
(Dynamics of the transformed system)
Suppose that $\sigma(A)\cap i\mathbb R=\emptyset$,  where $\sigma(A)$ is
the spectrum of the operator $A=\overline{f}'(x_0)$. Then
\begin{enumerate}
\item[$1)$] For the transformed equation $(\ref{eq3.21})$,
there exists a unique function
$\tilde{u}_{\varepsilon}\in C(\Omega,B[0,\beta])$ such that
\begin{equation} \label{eq3.24}
\tilde{u}_{\varepsilon}(\omega \tau)=\psi_{\varepsilon}(\tau,
\tilde{u}_{\varepsilon}(\omega),\omega)
\end{equation}
for all $\tau\in\mathbb R_{+}$ and $\omega\in\Omega$, where
$\psi_{\varepsilon}(\cdot,x,\omega)$ is a  unique solution of equation
(\ref{eq3.21}) which initial condition $\psi_{\varepsilon}(0,x,\omega)=x$;
\item[$2)$]
\begin{equation} \label{eq3.25}
\lim\limits_{\varepsilon\to 0}\max\limits_{\omega\in\Omega}
\Vert\tilde{u}_{\varepsilon}(\omega)\Vert=0.
\end{equation}
\end{enumerate}
\end{theorem}

\begin{proof} This statement follows from Theorem \ref{t2.6}.
\end{proof}

\begin{theorem}\label{t3.4}
(Analog of the second  Bogolyubov's theorem)
Assume that the conditions
$(i) - (iii)$ and $(\ref{eq3.6})$ are fulfilled and
$\sigma(A)\cap i\mathbb R=\emptyset$,  where $\sigma(A)$ is
the spectrum of the operator $A=\overline{f}'(x_0)$. Then for sufficiently
small $r_0>0$, there is $\varepsilon'_0$ with
 $0<\varepsilon'_0\le\varepsilon_0$ such that for
 $0<\varepsilon<\varepsilon'_0$, there exists a unique function
 $u_{\varepsilon}\in C(\Omega,B[x_0,r])$ such that
\begin{equation} \label{eq3.26}
{u}_{\varepsilon}(\omega t)=\psi_{\varepsilon}(t,{u}_{\varepsilon}(\omega),\omega)
\end{equation}
for all $t\in{\mathbb R}_+$ and $\omega\in\Omega$ and
\begin{equation} \label{eq3.27}
\lim\limits_{\varepsilon\to 0}\max\limits_{\omega\in\Omega}
\Vert{u}_{\varepsilon}(\omega)-x_0\Vert=0,
\end{equation}
where
$\psi_{\varepsilon}(\cdot,x,\omega)$ is the unique solution of the
nonautonomous equation
$(\ref{eq3.1})$ with  initial condition $\psi_{\varepsilon}(0,x,\omega)=x$,
and $x_0$ is a stationary solution of  the averaged equation $(\ref{eq3.3})$.
Note that the graph of ${u}_{\varepsilon}$ is an invariant integral manifold
for the nonautonomous equation $(\ref{eq3.1})$.
\end{theorem}

\begin{proof}  Under the conditions of the theorem and  in virtue of Theorem
\ref{t3.3} for the equation (\ref{eq3.21}),  there exists a unique function
$\tilde{u}_{\varepsilon}\in C(\Omega,B[0,r_0])$ with the properties
(\ref{eq3.24}) and (\ref{eq3.25}). Denote by
\begin{equation} \label{eq3.28}
{u}_{\varepsilon}(\omega)=x_0+\tilde{u}_{\varepsilon}(\omega)
-\varepsilon v(\omega,\tilde{u}_{\varepsilon}(\omega),\omega).
\end{equation}
Then from the equalities (\ref{eq3.14}), (\ref{eq3.15}),
(\ref{eq3.23}) and (\ref{eq3.27}),  we obtain the equality (\ref{eq3.28})
and the continuity of $u_{\varepsilon}:\Omega\to E$.
Consequently, we have,
 $u_{\varepsilon}\in C(\Omega,B[x_0,r])$  for sufficient small $\varepsilon>0$.
 The equality (\ref{eq3.26}) follows from the equalities
(\ref{eq3.24}) and (\ref{eq3.28}). The theorem is thus proved.
\end{proof}

\section{Almost periodic and recurrent solutions }

Let $ \mathbb T =\mathbb R $ or $\mathbb R_{+},\ (X,\mathbb T,\pi)$
be a dynamical system, $x\in X$, $\tau,\,\varepsilon
\in \mathbb T$, $\tau>0$, $\varepsilon>0$. We denote $\pi (x, t)$ by a short-hand
notation $xt$.

The point $x$ is called a stationary point if $xt=x$ for all $t\in \mathbb T$.
The point $x$ is called $\tau$-periodic if $x\tau=x$.

The number $\tau$ is called $\varepsilon$-shift ($\varepsilon$-almost period)
of a point $x$ if $\rho(x\tau,x)<\varepsilon$ ($\rho(x(t+\tau),xt)<\varepsilon$ for
all $t\in \mathbb T$).

The point $x$ is called  almost recurrent (almost periodic) if for any
$\varepsilon>0$ there exists positive number $l$ such that on every segment
of length $l$ can be found a $\varepsilon$-shift ($\varepsilon$-almost period)
of the point $x$.

A point $x$ is called recurrent if it is almost recurrent and the set
$H(x)=\overline{\{xt\mid t\in \mathbb T\}}$ is compact.

Denote by ${\frak M}_x=\{\{t_n\}\mid \{xt_n\} \,\,{\mbox is\,\, convergent}\}$.

\begin{theorem}\label{t4.1} (\cite{Sch1}, \cite{Sch2})
Let $(X,{\mathbb T}_1,\pi)$ and $(Y,{\mathbb T}_2,\sigma)$ be
dynamical systems with $ {\mathbb T}_1\subset{\mathbb T}_2$. Assume that
 $h:  X\to Y$ is a homomorphism from $(X,{\mathbb T}_1,\pi)$ onto
$(Y,{\mathbb T}_2,\sigma)$. If the point $x\in X$ is stationary
$(\tau$-periodic, quasi-periodic, almost periodic, recurrent$)$,
then the point $h(x)=y$ is also stationary $(\tau$-periodic,
quasi-periodic, almost periodic, recurrent$)$ and
${\frak M}_x\subset{\frak M}_y$.
 \end{theorem}

Consider the following nonautonomous equation in Banach space $E$
\begin{equation} \label{eq4.1}
w'=f(\omega t,w)
\end{equation}
where $f\in C(\Omega\times E,E)$. Suppose that the function
$f$ is regular,  i.e.,
 for all $w\in E$ and $\omega\in\Omega$, the equation (\ref{eq4.1}) admits a
 unique  solution $\varphi(t,w,\omega)$ defined on ${\mathbb R}_+$ with the
 initial condition $\varphi(0,w,\omega)=w$ and the mapping
 $\varphi:{\mathbb R}_+\times E\times \Omega\to E$ is continuous.

 It is well-known (see, for example, \cite{Sel}) that the
 mapping $\varphi$ satisfies the following conditions:

 a. $\varphi(0,w,\omega)=w$ for all $w\in E$ and $\omega\in\Omega$;

 b. $\varphi(t+\tau,w,\omega)=\varphi(t,\varphi(\tau,w,\omega),\omega \tau)$
 for all  $t,\tau\in{\mathbb T}_1 $, $w\in E$ and $\omega\in\Omega$.

 The solution $\varphi(t,w,\omega)$ of the equation (\ref{eq4.1}) is said
 to be
 stationary ($\tau$-periodic,  almost periodic, recurrent) if the point
 $x:=(w,\omega)\in X:=E\times \Omega$ is stationary ($\tau$-periodic,  almost
 periodic, recurrent)
 point of the skew-product   dynamical system $(X,{\mathbb R}_+,\pi)$, where
 $\pi=(\varphi,\sigma)$, i.e. $\pi(t,(w,\omega))=(\varphi(t,w,\omega),\omega t)$
 for all $t\in{\mathbb R}_+$ and $(w,\omega)\in E\times \Omega$.

\begin{lemma}\label{l4.2}
Suppose that $u\in C(\Omega,E)$  satisfies  the condition
\begin{equation} \label{eq4.2}
u(\omega t)=\varphi(t,u(\omega),\omega)
\end{equation}
for all $t\in{\mathbb R}$ and $\omega\in\Omega$. Then the mapping
$h:\Omega\to X$ defined by
\begin{equation} \label{eq4.3}
h(\omega)=(u(\omega),\omega)
\end{equation}
for all $\omega\in\Omega$ is a homomorphism  from
$(\Omega,{\mathbb R},\sigma)$
onto
$(X,{\mathbb R}_+,\pi)$.
\end{lemma}

\begin{proof} This assertion follows from the equalities (\ref{eq4.2}) and
(\ref{eq4.3})
\end{proof}

\begin{remark} The function $u\in C(\Omega,E)$ with the property $(\ref{eq4.2})$
is called continuous invariant section (or integral manifold) for
non autonomous system $(\ref{eq4.1})$.
\end{remark}

\begin{theorem}\label{t4.3} If the function $u\in C(\Omega,E)$
satisfies the condition $(\ref{eq4.2})$ and the point $\omega\in\Omega$
is stationary $(\tau$-periodic,  almost periodic, recurrent$)$, then
the solution  $\varphi(t,u(\omega),\omega)$ of the equation $(\ref{eq4.1})$
 also will be stationary $(\tau$-periodic,  almost periodic, recurrent$)$.
 \end{theorem}

\begin{proof} This statement follows from Theorem \ref{t4.1} and Lemma
\ref{l4.2}.
\end{proof}

\begin{example}\label{ex4.4}{\rm Consider the equation
\begin{equation} \label{eq4.4}
u'=f(t,u)
\end{equation}
where $f\in C(\mathbb R\times E,E)$; here $ C(\mathbb R\times E,E)$ is the space of all
continuous function $\mathbb R\times E\to E)$ equipped with compact-open
topology. Along with the equation $(\ref{eq4.4})$,  we will consider
the $H$-class
of equation $(\ref{eq4.4})$
\begin{equation} \label{eq4.5}
u'=g(t,u) \quad (g\in H(f)),
\end{equation}
where $H(f)=\overline{\{f_{\tau}\mid \tau\in\mathbb R\}}$ and the over bar
denotes
the closure in $C(\mathbb R\times E,E)$ and $f_{\tau}(t,u)=f(t+\tau,u)$ for
all $t\in{\mathbb R}$ and $u\in E$. Denote by
$(C(\mathbb R\times E,E),\mathbb R,\sigma)$ the Bebutov's dynamical system
(see, for example, \cite{Sch1}, \cite{Sch2}, \cite{Sel}). Here
$\sigma(t,g)=g_t$ for all $t\in{\mathbb R}$ and $g\in C(\mathbb R\times E,E)$.
}
\end{example}

The function $f\in C(\mathbb R\times E,E)$ is called regular (see \cite{Sel})
if for all $u\in  E $ and $g\in H(f)$ the equation (\ref{eq4.5}) admits a
unique solution $\varphi(t,u,g)$ defined on ${\mathbb R}_+$ with the initial
condition $\varphi(0,u,g)=u$.

Let $\Omega$ be the hull $H(f)$ of a given regular function
$f\in C(\mathbb R\times E,E)$ and denote the restriction
of $(C(\mathbb R\times E,E),\mathbb R,\sigma)$ on $\Omega$ by
$(\Omega,\mathbb R,\sigma)$. Let $F:\Omega\times E\to E$ be a continuous
mapping defined by $F(g,u)=g(0,u)$ for $g\in \Omega$ and $u\in E$. Then the
equation (\ref{eq4.5}) can be written in such form:
\begin{equation} \label{eq4.6}
u'=F(\omega t,u),
\end{equation}
where $\omega=g$ and $\omega t=g_t$.
\begin{lemma}\label{l4.5}
The following two conditions are equivalent.
\begin{enumerate}
\item[$1)$] There exists a limit
\begin{equation} \label{eq4.7}
 f_0(x)=\lim\limits_{T\to+\infty}\frac{1}{T}\int\limits_t^{t+T}f(s,x)ds
\end{equation}
 uniformly w.r.t. $t\in \mathbb R$ and $x$ on every compact set $K\subset E$.
\item[$2)$] There exists a limit
\begin{equation} \label{eq4.8}
 f_0(x)=\lim\limits_{T\to+\infty}\frac{1}{T}\int\limits_0^{T}g(s,x)ds
\end{equation}
 uniformly w.r.t. $g\in H(f)$ and $x$ on every compact set $K\subset E$.
\end{enumerate}
\end{lemma}

\begin{proof}  The equality (\ref{eq4.7}) follows from (\ref{eq4.8}) because
$f_t\in H(f)$ for all $t\in \mathbb R$ and
$$\lim\limits_{T\to+\infty}\frac{1}{T}\int\limits_0^{T}f_t(s,x)ds=
 \lim\limits_{T\to+\infty}\frac{1}{T}\int\limits_t^{t+T}f(s,x)ds.$$

Let $g\in H(f)$, then there exists a sequence $\{t_n\}\subset \mathbb R$ such
that $g=\lim\limits_{n\to+\infty}f_{t_n}$. From the equality (\ref{eq4.7}), it
follows that for all $\varepsilon>0$ and compact set $K\subset E$, there
exists $L(\varepsilon, K)>0$  such that
\begin{equation} \label{eq4.9}
\Vert \frac{1}{T}\int\limits_0^{T}f_{t_n}(s,x)ds-f_0(x)\Vert<\varepsilon
\end{equation}
for all $n\in \mathbb N$ and $T\ge L(\varepsilon,K)$. Passing to the
limit in the
equality (\ref{eq4.9}) as $n\to+\infty$ we obtain the equality (\ref{eq4.8}).
\end{proof}

\begin{theorem}\label{t4.6}
(Recurrent solutions)
Suppose that the following conditions are fulfilled:
\begin{enumerate}
\item[$1)$] $f\in C({\mathbb R}\times E,E)$ and there exist $x_0\in E$ and
$r>0$ such that the function $f$ is bounded on ${\mathbb R}\times B[x_0,r]$,
i.e., there exists positive number $M$  such that
\begin{equation} \label{eq4.10}
\Vert f(t,x)\Vert\le M
\end{equation}
for all $t\in{\mathbb R}$ and $x\in B[x_0,r]$.
\item[$2)$] The functions $f\in C({\mathbb R}\times E,E)$ and
$f_0\in C(E,E)$  are twice continuously differentiable w.r.t. variable
$x\in B[x_0,r]$. Moreover,  the function $f'_x(t,x)$ is
bounded on ${\mathbb R}\times B[x_0,r]$, and $f'_0(x)$
is bounded on $B[x_0,r]$.
\item[$3)$] The equality (\ref{eq4.7}) can be twice differentiated, i.e.
the following equalities
\begin{equation} \label{eq4.11}
 {f}_0^{(k)}(x)=
 \lim\limits_{T\to +\infty}\frac{1}{T}\int\limits_t^{t+T}f_x^{(k)}(s,x)ds
 \quad (k=1,2)
\end{equation}
 take place, uniformly w.r.t. $t\in\mathbb R$ and $x\in B[x_0,r]$.
\item[$4)$] $f_0(x_0)=0$ and $\sigma(A)\cap i{\mathbb R}=\emptyset$, where
$A=f'_0(x_0)$ and   $\sigma(A)$ is  the spectrum of operator $A$.
\item[$5)$] The function $f\in C({\mathbb R}\times E,E)$
is stationary $(\tau$-periodic,  almost periodic, recurrent$)$
w.r.t. $t\in{\mathbb R}$, and  uniformly w.r.t. to $x$ on every
compact subset $K \subset E$.
\end{enumerate}

Then for sufficiently  small $r_0>0$, there exists $0<\varepsilon'_0\le\varepsilon_0$
such that for $0<\varepsilon<\varepsilon'_0$ the equation
\begin{equation}\label{eq4.12}
x'=\varepsilon f(t,x)
\end{equation}
admits a unique stationary $(\tau$-periodic, almost periodic, recurrent$)$
solution $\varphi^{\varepsilon}(t)$ with the following properties:
\begin{enumerate}
\item[$a)$] $\Vert \varphi^{\varepsilon}(t)-x_0\Vert\le r_0$ for all
$t\in{\mathbb R}$.
\item[$b)$] $\lim\limits_{\varepsilon\to 0}\sup\limits_{t\in{\mathbb R}}
\Vert \varphi^{\varepsilon}(t)-x_0\Vert=0$.
\item[$c)$] ${\frak M}_f\subset {\frak M}_{\varphi^{\varepsilon}}$, where
${\frak M}_f=\{\{t_n\}\mid \{f_{t_n}\}\,\,{\mbox is\,\, convergent\,\, on\,\,}
C({\mathbb R}\times E,E)\}$ and
${\frak M}_{\varphi^{\varepsilon}}=\{\{t_n\}\mid
\{\varphi^{\varepsilon}_{t_n}\}\,\,{\mbox is\,\, convergent\,\, on\,\,}
C({\mathbb R},E)\}$.
\end{enumerate}
\end{theorem}

\begin{proof}  Note that   the conditions 1) -- 4) imply
conditions
$(i)$ -- $(iii)$.  So the proof of the theorem
follows from Theorems \ref{t3.4}, \ref{t4.1}, \ref{t4.3} and
Lemma \ref{l4.5}.
\end{proof}

\begin{remark}
Note that  Theorem $\ref{t4.6}$ is also true for the equation
$(\ref{eq4.12})$ with non recurrent function $f$. For example, if
$f$ is pseudo recurrent  (\cite{Sch1,Sch2}), i.e., if $H(f)$
is compact and every function $g\in H(f)$ is stable in the sense
of Poisson. In this case we can affirm that the solution
$\varphi^{\varepsilon}$ will be also pseudo recurrent.
See also Section 6 later in this paper.
\end{remark}

\section{Invariant torus and quasi-periodic solutions}

Let ${\mathcal T}^m $ be an $m$-dimensional torus. We consider a
  nonautonomous dynamical system in Banach space $E$,
with a driving system defined on the torus  ${\mathcal T}^m $:
\begin{equation} \label{eq5.1}
\left\{\begin{array}{rl}
x'=A(\omega)x+f(\omega)+F(\omega,x)\\
\omega'=\Phi(\omega),
\end{array}
\right.
\end{equation}
where $\Phi\in C({\mathcal T}^m,T{\mathcal T}^m)$, $T{\mathcal T}^m$ is a
tangent space of the torus ${\mathcal T}^m$, $f\in C({\mathcal T}^m,E)$,
$A\in C(\Omega, L(E) )$ and $F\in C({\mathcal T}^m\times E,E)$.

We suppose that the second equation of the system (\ref{eq5.1})
generates an autonomous dynamical system
$({\mathcal T}^m,{\mathbb R},\sigma)$
on the torus ${\mathcal T}^m$ and the equation
\begin{equation} \label{eq5.2}
x'=A(\omega t)x+f(\omega t)+F(\omega t,x)
\end{equation}
admits a unique solution $\varphi(t,x,\omega)$ defined on ${\mathbb R}_+$ and
satisfying the initial condition $\varphi(0,x,\omega)=x$.

A function $\gamma\in C({\mathcal T}^m,E)$ is called \cite{Sam}
an $m$-dimensional invariant torus of equation (\ref{eq5.2})
(or system (\ref{eq5.1})) if
\begin{equation} \label{eq5.3}
\gamma(\omega t)=\varphi(t,\gamma(\omega),\omega)
\end{equation}
for all $t\in{\mathbb R}_+$ and $\omega\in {\mathcal T}^m$.

Applying the results from sections 2 -- 4 ,
 we  have the following tests
of existence of the invariant torus for   equation (\ref{eq5.2}).

\begin{theorem}\label{t5.1}
(Invariant torus)
Suppose that the  equation (\ref{eq2.1}) is hyperbolic
and there exist positives numbers $0<L<L_0:=\frac{\nu}{2N}$ and
$0<r<r_0:=\nu_0(\frac{\nu}{2N}-L_0)^{-1}$ such that the function
$F\in C({\mathcal T}^m\times E,E)$ satisfies the condition (\ref{eq2.6}).
Then the equation (\ref{eq5.2}) admits an $m$-dimensional invariant torus.
\end{theorem}

\begin{theorem}\label{t5.2}
(Invariant torus for perturbed system)
Suppose that there exist positives numbers $r$ and
$L$ such that the  condition (\ref{eq2.13}) is fulfilled.
Then for sufficient small $\varepsilon\ge0$ there exists
an $m$-dimensional invariant torus $u_{\varepsilon}$ for the perturbed
equation
\begin{equation} \label{eq5.4}
x'=A(\omega t)x+f(\omega t)+\varepsilon F(\omega t,x)\quad
(\omega\in {\mathcal T}^m)
\end{equation}
and
\begin{displaymath}
\lim\limits_{\varepsilon\to 0}\max\limits_{\omega\in\Omega}
\Vert u_{\varepsilon}(\omega)-u_0(\omega)\Vert=0.
\end{displaymath}
\end{theorem}

\begin{theorem}\label{t5.3}
(Unique invariant torus)
Let $\Omega={\mathcal T}^m$.  Assume  the conditions
of Theorem \ref{t3.4} are satisfied. Then, for a sufficient small $r_0>0$,
 there exists $0<\varepsilon_0'<\varepsilon_0$
 such that for all $0<\varepsilon<\varepsilon_0'$ there exists
a unique $m$-dimensional invariant torus $u_{\varepsilon}$ for the
equation (\ref{eq5.1}) and
\begin{displaymath}
\lim\limits_{\varepsilon\to 0}\max\limits_{\omega\in\Omega}
\Vert u_{\varepsilon}(\omega)-x_0\Vert=0.
\end{displaymath}
\end{theorem}

We have the following corollary for quasi-periodic
nonautonomous dynamical systems, i.e., the driving system defined on the
torus ${\mathcal T}^m $ is
quasi-periodic in time.

\begin{coro}\label{c5.4}
(Compact minimal invariant torus)
Suppose that the conditions of Theorem \ref{t5.1}
(respectively Theorem \ref{t5.2} or Theorem \ref{t5.3}) are
fulfilled and the dynamical system $({\mathcal T}^m,{\mathbb
R},\sigma)$ generated by the second equation of the system
(\ref{eq5.1}) is compact minimal and contains only quasi-periodic
motions, then the equation (\ref{eq5.2}) (respectively the
equation (\ref{eq5.4}) or the equation (\ref{eq3.1})) admits an
$m$-dimensional invariant torus $u_{\varepsilon}$ which is compact
minimal and contains only quasi-periodic motions.
\end{coro}

\section{Pseudo recurrent solutions}

An autonomous dynamical system $(\Omega, {\mathbb T},\sigma)$
 is said to be
pseudo recurrent   if the following conditions are fulfilled:
\begin{enumerate}
\item[a)]
$\Omega$  is compact;
\item[b)]
$(\Omega, {\mathbb T},\sigma)$ is transitive, i.e. there exists a point
$\omega_0\in\Omega$ such that
$\Omega=\overline{\{\omega_0t\mid t\in {\mathbb T}\}}$;
\item[c)]
every point $\omega\in\Omega$ is stable in the sense of Poisson, i.e.
\begin{displaymath}
{\frak N}_{\omega}=\{\{t_n\}\mid \omega t_n\to \omega\; {\mbox and}\;
\vert t_n\vert\to +\infty\}\ne\emptyset.
\end{displaymath}
\end{enumerate}

\begin{lemma}\label{l6.1}
Let $<(X,{\mathbb T}_1,\pi),(\Omega,{\mathbb T}_2,\sigma),h >$
be a nonautonomous dynamical system and the following conditions
are fulfilled: \been
\item[$1)$] $(\Omega,{\mathbb T}_2,\sigma)$ is pseudo recurrent;
\item[$2)$] $\gamma\in C(\Omega,X)$ is an invariant section of the
homomorphism $h:X\to\Omega$.
\enen

Then the autonomous dynamical system
$(\gamma(\Omega),{\mathbb T}_2,\pi)$ is  pseudo recurrent.
\end{lemma}

\begin{proof} It is evident that the space $\gamma(\Omega)$ is compact,
because $\Omega$ is compact and $\gamma\in C(\Omega,X)$.
We note that on the space $\gamma(\Omega)$, by the
homomorphism $\gamma:\Omega\to\gamma(\Omega)$, we have    a
 dynamical
system $(\gamma(\Omega),{\mathbb T}_2,\widehat{\pi})$, namely
$\widehat{\pi}^t\gamma(\omega):=\gamma(\omega t)$ for all $t\in{\mathbb T}_2$
and $\omega\in\Omega$, then $\widehat{\pi}^t\gamma(\omega)={\pi}^t\gamma(\omega)$
for all $t\in{\mathbb T}_1\subseteq{\mathbb T}_2$ and $\omega\in\Omega$.
Now we will show that
$\gamma(\Omega)=\overline{\{\gamma(\omega_0) t\mid
t\in{\mathbb T}_2\}}$.
In fact, let $x\in\gamma(\Omega)$. Then there exists
a unique point $\omega\in\Omega$ such that $x=\gamma(\omega)$. Let
$\{t_n\}\subset {\mathbb T}_2$ be a sequence such that $\omega_0t_n\to \omega$.
Then $x=\gamma(\omega)= \lim\limits_{n\to+\infty}\gamma(\omega_0t_n)=
\lim\limits_{n\to+\infty}\gamma(\omega)t_n)$ and, consequently,
$\gamma(\Omega)\subset\overline{\{\gamma(\omega_0) t\mid t\in{\mathbb T}_2\}}$.
The inverse inclusion is trivial. Hence,
$\gamma(\Omega)=\overline{\{\gamma(\omega_0) t\mid
t\in{\mathbb T}_2\}}$.
To finish the proof of the lemma it is sufficient to note that
${\frak N}_{\omega}\subseteq{\frak N}_{\gamma(\omega} $
for every point $\omega\in\Omega$ and, consequently, every point
$\gamma(\omega)$ is Poisson stable. The lemma is proved.
\end{proof}

  Lemma \ref{l6.1} implies  that   the conditions of Theorem
\ref{t5.1} (respectively Th. \ref{t5.2} or Th. \ref{t5.3})  are satisfied.
Therefore, we have the following result.

\begin{theorem}\label{soudo}
(Pseudo recurrent integral manifold)
Assume the driving dynamical system  $(\Omega, {\mathbb T},\sigma)$
 is  pseudo recurrent , and assume the conditions in Lemma \ref{l6.1} are satisfied. Then
the
equation (\ref{eq5.2}) (respectively, the equation (\ref{eq5.4}) or
the equation (\ref{eq3.1})) admits a pseudo recurrent integral
manifold.
\end{theorem}

\bigskip

{\bf Acknowledgment:} The research described in this publication
was made possible in part by Award No. MM1-3016 of the Moldovan
Research and Development Association (MRDA) and  the U.S. Civilian
Research \& Development Foundation for the Independent States of
the Former Soviet Union (CRDF). This paper was written while the
first author was visiting   Illinois Institute of Technology
(Department of Applied Mathematics) in Spring of 2002. He would
like to thank people in that institution for their very kind
hospitality.

\end{document}